\documentclass[11pt]{article}
\usepackage{amsfonts,amsmath,amssymb,amsthm, amscd}
\usepackage[dvips]{epsfig,graphics}

\numberwithin{equation}{section}

\newtheorem{theorem}{Theorem}[section]
\newtheorem{prop}[theorem]{Proposition}

\newtheorem{cor}[theorem]{Corollary}

\theoremstyle{definition}
\newtheorem{definition}[theorem]{Definition}
\newtheorem{example}[theorem]{Example}
\newtheorem{remark}[theorem]{Remark}

\def\<{{\langle}}
\def\>{{\rangle}}
\def\G{{\Gamma}}

\def\Z{\mathbb Z}
\def\S{{\Sigma}}

\def\sx{\sigma_x}
\def\G{{\cal G}}

\begin{document}

\title{Nonfibered knots and representation shifts}

\author{Daniel S. Silver \and Susan G. Williams\thanks{Both authors partially supported by NSF grant
DMS-0304971.} \\ {\em
{\small Department of Mathematics and Statistics, University of South Alabama}}}

\maketitle 

\begin{abstract}\noindent A conjecture of \cite{swTAMS} states that a knot is nonfibered if and only if its infinite cyclic cover has uncountably many finite covers. We prove the conjecture for a class of knots that includes all knots of genus 1, using techniques from symbolic dynamics.
\end{abstract}

 \begin{footnote}{Mathematics Subject Classification:  
Primary 57M25; secondary 37B40.}\end {footnote}

\section{Introduction}
\label{intro} 

Let $G$ be a finitely presented group with epimorphism $\chi: G \to \Z$. The  kernel $K$ of $\chi$ need not be finitely generated. However, $K$ is finitely presented as a ${\mathbb Z}$-operator group \cite{swIJM}. In \cite{swJPAA} \cite{swTAMS} the authors exploited this structure to show that the representations of $K$ into a fixed finite  group $\S$ form a {\it shift of finite type}, a simple dynamical system described by a finite directed graph. We call this dynamical system the {\it representation shift} of $K$ in $\S$. When $G$ is a knot or link group, representation shifts inform us about the algebraic topology of finite covering spaces from a purely dynamical perspective. 

We review basic definitions of representation shifts and give a partial solution to Conjecture 4.4 of \cite{swTAMS}. The complete solution would characterize nonfibered knots as knots with complicated representation shifts, where complexity is measured by topological entropy.

Part of this work was done during the M.M. Postnikov Memorial Conference, Bedlewo, Poland in June 2007. We are grateful to the conference organizers for giving us the opportunity to participate. We thank Stefan Friedl for stimulating and helpful discussions. 

\section{Review of Representation Shifts}
\label{rep}

An {\it augmented group system} \cite{sMA} is a triple $\G=(G, \chi, x)$ consisting of a finitely presented group $G$, epimorphism $\chi: G \to {\mathbb Z}$ and distinguished element $x \in G$ such that $\chi(x)=1$. Two such systems $\G_i=(G_i, \chi_i, x_i), \ i=1,2$ are {\it equivalent} (and regarded as the same) if there exists an isomorphism $f:G_1\to  G_2$ such that $f(x_1)=x_2$ and $\chi_1 = \chi_2\circ f$. 

\begin{example} \label{ags} An augmented group system is associated to an oriented knot $k\subset {\mathbb S}^3$ in a canonical manner. Let $G=\pi_1({\mathbb S}^3\setminus k, p)$, where the base point $p$ is contained on the boundary $\partial N(k)$ of a tubular neighborhood $N(k)= {\mathbb S}^1 \times {\mathbb D}^2$ of $k$. Let $x$ be the homotopy class of a meridian $m\subset \partial N(k)$, with orientation acquired from $k$. Finally, let $\chi: G \to {\mathbb Z}$ be the abelianization homomorphism that sends $x$ to $1$.  It follows from the uniqueness of tubular neighborhoods that $\G=(G, \chi, x)$ is well defined. \end{example}

We denote the kernel of $\chi$ by $K$. Given any finite group $\S$, we consider the space 
${\rm Hom}(K, \S)$ of representations $\rho: K \to \S$. The basis for its topology is given by the sets
$${\cal N}_{a_1, \ldots, a_s}(\rho)=\{\rho' \mid \rho'(a_i)=\rho(a_i),\ i=1, \ldots, s\},$$
where $a_1, \ldots, a_s$ varies over all finite collections of elements of $K$. The topology is the compact-open topology where $K$ and $\S$ are discrete spaces. Roughly speaking, representations are close in ${\rm Hom}(K, \S)$ if they agree on large finitely generated subgroups of $K$. The distinguished element $x$ induces a self-map $\sx$ of ${\rm Hom}(K, \S)$ defined by 
$$\sx \rho(a) = \rho(x^{-1}ax)\quad \forall a \in K.$$
It is easily seen that $\sx$ is a homeomorphism. 

The {\it representation shift} associated to $\G=(G, \chi, x)$ and $\S$ is the pair $({\rm Hom}(K, \S), \sx).$  We denote it by $\Phi_\S(\G)$. 
It is a dynamical system  well defined up to topological conjugacy  \cite{swIJM}. More precisely, if $\G_i,\ i=1,2$, are equivalent augmented group systems, then there exists a homeomorphism $F$ of the underlying spaces of $\Phi_\S(\G_i)$ such that $F\circ \sigma_{x_1}=\sigma_{x_2}\circ F$. 

The representation shift $\Phi_\S(\G)$ is an example of a {\it shift of finite type}, a special type of expansive $0$-dimensional dynamical system, one that can be described by a finite directed graph. (See \cite{lm}.) We use combinatorial group theory to construct such a graph for a representation shift.  

Given an augmented group system $\G=
(G, \chi, x)$, we can describe $G$ as an HNN extension $\<x, B\mid x^{-1}ax=\phi(a),\ \forall a \in U\>$, where $B$ is a finitely generated subgroup of $K$, and $U$, $V$ are isomorphic finitely generated subgroups of $B$ with isomorphism $\phi: U \to V$ (see \cite{ls}). The subgroup $B$ is an
{\it HNN base}. One can choose $B$ so that it contains any prescribed finite subset of $K$ (see \cite{sPAMS}).

\begin{example} \label{hnn} Let $\G=(G, \chi, x)$ be an augmented group system associated to a knot, as in Example \ref{ags}.  An HNN decomposition for $G$ can be obtained in a natural way. Begin with a $\pi_1$-incompressible Seifert surface for $k$ meeting the exterior $E(k)={\mathbb S}^3\setminus {\rm int}\ N(k)$ in a connected surface $S$. Let $(W; S_0, S_1)$ be the resulting cobordism, with boundary comprising two copies $S_0, S_1$ of $S$ joined by an annulus $\partial S \times I$. Let $B= \pi_1(W, p)$, where the basepoint $p$ lies on the boundary of $S_0$. Let $U= \pi_1(S_0, p)$. The meridian $m$ appears as a path from $p \in S_0$ to a point $p_1 \in S_1$. Use the path to regard $\pi_1(S_1, p_1)$ as a subgroup $V$ of $B$. Clearly $G$ is described as $(B; U, V, \phi)$, where $\phi$ is induced by the gluing of $S_0$ to $S_1$ when recovering the exterior $E(k)$. \end{example}

Conjugation by $x$ induces an automorphism of $K$. Let $B_j = x^{-j}Bx^j$, $U_j = x^{-j}Ux^j$ and $V_j = x^{-j}V x^j$, for $j \in {\mathbb Z}$. Then $K$ is described as an infinite amalgamated free product
$$K = \<B_j\mid V_j = U_{j+1},\ \forall j \in {\mathbb Z}\>.$$

The vertex set of the graph $\Gamma$ consists of all representations $\rho_0: U \to \S$, a finite set since $U$ is finitely generated. If $\bar \rho_0$ is a representation from $B$ to $\S$, then we draw a directed edge labeled $\bar \rho_0$ from the vertex $\rho_0 =\bar \rho_0|_U$ to the vertex $\rho_0' = \bar \rho_0|_V\circ \phi$. ($\Gamma$ may have parallel edges.) Consider a bi-infinite path in $\Gamma$ given by an edge sequence
$$\cdots\ \bar \rho_{-2}\ \bar\rho_{-1}\ \bar\rho_0\ \bar\rho_1\ \bar\rho_2\ \cdots$$
The representations $B_j \to \S$ given by $a \mapsto \bar\rho_j(x^jax^{-j})$ have a unique common extension $\rho: K \to \S$. Conversely, any representation $\rho: K \to \S$ arises from such a path, and uniquely.
Thus bi-infinite paths of the graph $\Gamma$ correspond bijectively to elements of ${\rm Hom}(K, \S)$. The map $\sx$ acts as the left coordinate shift on the sequence of edges. 

We may ``prune" $\Gamma$ by removing any vertex or edge that is not contained in a bi-infinite path.  The resulting graph has finitely many bi-infinite paths iff it consists of a collection of disjoint cycles. It contains uncountably many bi-infinite paths iff it contains two cycles with at least one common vertex.  

A representation $\rho \in \Phi_\S(\G)$ has {\it period} $r$ if $\sx^r(\rho)= \rho$. Such representations correspond to closed paths in $\Gamma$ with 
length dividing $r$. The set of representations with period $r$ is denoted by ${\rm Fix}(\sx^r)$. If $M_r$ is the $r$-fold cyclic cover of ${\mathbb S}^3$ branched over a knot $k$, then ${\rm Fix}(\sx^r)$ is in natural bijective correspondence with ${\rm Hom}(\pi_1 M_r, \S)$ \cite{swTAMS}. This correspondence connects dynamical properties of  the representation shift with topological properties of $k$.

Topological entropy is one measure of  complexity for a dynamical system. 
For a shift of finite type, it can be computed as the log of the spectral radius of the adjacency matrix $A$ of any directed graph that describes the shift. (Here $A_{i,j}$ is the number of edges from the $i$th  vertex to the $j$th.) Consequently, the topological entropy of $\Phi_\S(\G)$, denoted by $h_\S(\G)$, is the exponential growth 
rate of $|{\rm Hom}(\pi_1 M_r, \S)|$ (see \cite{swTAMS}). Notice that if $K$ is finitely generated, then $\Phi_\S(\G)$ is finite for all $\S$, and so in this  case $h_\S(\G)$=0.

Let $S_N$ denote the symmetric group on $\{1, \ldots, N\}$. It is well known that elements $\rho \in {\rm Hom}(K, S_N)$ correspond in a finite-to-one manner with subgroups $H\le K$ with index no greater than $N$. The correspondence is 
$$\rho \mapsto \{g\in K \mid \rho(g)(1)=1\}.$$
The preimage of a subgroup of index $N$ consists of $(N-1)!$ transitive representations. (A representation $\rho$ is  transitive if $\rho(K)$  operates transitively on $\{1, \ldots, N\}$.) Note that if $\Phi_{S_N}(\G)$ is uncountable, then $K$ contains uncountably many subgroups of some index no greater than $N$. Hence the infinite cyclic cover of  $k$ 
has uncountably many finite covers. 

We summarize the results of this section. Recall that any finite group embeds in a sufficiently large symmetric group.

\begin{prop} \label{summary} Let $k\subset {\mathbb S}^3$ be a knot with associated augmented group system $\G$. Then the following  statements are equivalent.
\item{(1)} The infinite cyclic cover of $k$ has  uncountably many finite covers. 
\item{(2)} The representation shift $\Phi_{\S}(\G)$ is uncountable, for some finite group $\S$. 
\item{(3)} The topological entropy $h_{\S}(\G)$ is positive, for some finite group $\S$.
\item{(4)} $\displaystyle \lim_{r\to \infty} \displaystyle{1\over r} \log |{\rm Hom}(\pi_1 M_r, \S)|$ is positive, for some finite group $\S$.

\end{prop}

\section{Nonfibered knots} \label{nonfibered} 
We recall that a knot $k \subset {\mathbb S}^3$ is fibered if its exterior
$E(k)= {\mathbb S}^3 \setminus  {\rm int}\ N(k)$ fibers over the circle. It is no loss of generality to assume that the fibration restricts to the standard projection
$\partial N(k) \simeq k \times {\mathbb S}^1 \to {\mathbb S}^1$. Hence
$E(k)$ is seen to be homeomorphic to a mapping torus $S\times I/F$, where $F: S \to S$ is a homeomorphism of a minimal-genus Seifert surface $S$ of $k$. 

If $k$ is fibered, then the commutator subgroup $G'$ of its group is finitely generated and free, isomorphic to $\pi_1 S$. Conversely, a theorem of J. Stallings \cite{stallings} implies that if $k$ is a knot such that $G'$ is finitely generated, then in fact $G'$ is free and $k$ is fibered.

If $k$ is fibered and $\G$ is its associated augmented group system, then for any finite group $\S$, the representation shift $\Phi_\S(\G)$ is finite. 
Its order is $|\S|^{2g}$, where $g$ is the genus of $k$ (equal to the genus of its fiber).  The trefoil and figure-eight knots are the only fibered knots of genus 1. 

Conjecture 4.4 of \cite{swTAMS} proposes a characterization of nonfibered knots. It states that $k$ is nonfibered iff the entropy $h_\S(\G)$ is positive for some finite group $\S$. 

\begin{remark} \label{rem} (1) In terms of the HNN base $B$ described above, the condition that $k$ is not fibered is equivalent to the condition that $U$ is a proper subgroup of $B$. Lemma 2.3 (Substitution Lemma) of \cite{swJPAA} provides a strategy for showing that some
$\Phi_{S_n}(\G)$ is uncountable: Find a periodic element of $\Phi_{S_N}(\G)$ such that some symbol, say $N$, is fixed by every permutation in the image of $U$ but moved by some element of $\rho(K)$. Recall that periodic representations correspond to cycles in the graph $\Gamma$. By introducing a new symbol (enlarging $S_N$ to $S_{N+1}$), we can construct another periodic representation corresponding to a second cycle, one that branches from the first. Then $\Phi_{S_{N+1}}(\G)$ is uncountable. 

(2) For our strategy, it suffices to find any representation $\tilde \rho: G \to  \S$ such that $\rho(U)$ is a proper subgroup of $\rho(K)$.  For given such 
a representation, and letting $\rho: K \to \S$ be the restriction, we enumerate the cosets of $\rho(U)$ in $\rho(K)$, say $1,\ldots, N$ ($N>1$). In a natural way, $\rho$ determines an element of  $\Phi_{S_N}(\G)$:  $a \in K$ is  sent to the transitive permutation of cosets given by right multiplication by $\rho(a)$. Note that if $a\in U$, then such a permutation fixes the symbol corresponding to $\rho(U)$. Finally, we note that if $r$ is the order of $\tilde \rho(x)$ in $\S$, then $\sx^r\rho = \rho$, since
$(\sx^r\rho)(a) = \rho(x^{-r}ax^r) = \tilde\rho(x^{-1})^r\rho(a)\tilde\rho(x)^r = \rho(a)$, for all $a\in K$. 
\end{remark}

The representation $\tilde \rho$ in the Remark \ref{rem} (2) ``separates" the subgroup 
$U$ from some element $a\in K$.

In general, a subgroup $U$ of a group $G$ is {\it separable} if for any element $a \in G\setminus U$, there exists a finite-index subgroup of $G$ that contains $U$ but not $a$. Equivalently, there exists a finite representation $\tilde\rho:  G \to \S$ such that $\tilde \rho(a) \notin \tilde\rho(U)$. The strategy outlined in Remark \ref{rem}(2) requires only that $U$  can be separated from {\sl some} element of $K\setminus U$.

\begin{definition} An 
element $a \in G\setminus U$ is {\it separable from} $U$ if there exists a subgroup $H$ of finite index in $G$ containing $U$ but not $a$.  \end{definition}

Question 15 of \cite{thurston} asks if any finitely generated subgroup of a finitely-generated Kleinian group is separable. An affirmative answer would establish Theorem \ref{main} for all hyperbolic knots.
Although Thurston's question remains open, a result  of D. Long and G. Niblo \cite{ln} enables us to apply our strategy in the case of genus-1 knots (see also remarks that follow). 

The theorem of Long and Niblo has been used by S. Friedl and S. Vidussi in \cite{fv} to show that twisted Alexander polynomials corresponding to finite representations decide if a genus-1 knot is fibered. 

\begin{theorem}\label{ln} 
{\rm (D. Long and G. Niblo \cite{ln})} Let $M$ be an orientable Haken $3$-manifold. If $i: T \hookrightarrow M$ is an incompressibly embedded torus, then 
$i_*(\pi_1 T)$ is separable in $\pi_1 M$. \end{theorem}

\begin{theorem} \label{main} Let $k$ be a knot of genus 1. Then $k$ is nonfibered iff the conclusions of Proposition \ref{summary} hold. 
\end{theorem}

\begin{proof} One implication of the theorem is clear: if the conclusion of Proposition \ref{summary} holds, then $k$ is nonfibered. 

Assume that $k$ is nonfibered. Consider the $3$-manifold $M$ obtained  by $0$-framed surgery on $k$; that is, by removing and replacing a tubular neighborhood $N(k)\equiv k \times {\mathbb D}^2$ in such a way that each disk $*\times {\mathbb D}^2$ bounds a longitude of $k$ . By results of \cite{gab2}, $M$ is irreducible. We denote the fundamental group of $M$ by $\hat G$. 

The addition of a meridianal disk converts a genus-1 Seifert surface $S$ for $k$ to a torus $\hat S$ in $M$.
Since $\hat S$ is dual to a nontrivial cohomology class and $M$ is irreducible, we see that $\hat S$ is incompressible. Note in particular that $M$ is Haken.

Obtain an HNN decomposition $(\hat B; \hat U, \hat V)$ for $\hat G$ much as we did for $G$, by 
splitting $M$ along $\hat S$.  Here $\hat U= \pi_1 \hat S$.
Since $k$ is not fibered, neither is  $M$ \cite{gab1}. Hence
$\hat U$ must be a proper subgroup of $\hat B$. Select an element  $\hat a\in \hat B\setminus \hat U$.   By Theorem \ref{ln} there exists a finite group $\S$ and homomorphism $\hat  \rho: \hat G \to \S$ such  that $\hat \rho(\hat a) \notin \hat U$. 

The group $\hat G$ is  a quotient of $G$. Let $p$ be the natural projection. 
Note that  $p(U)=\hat U$. Choose $a\in K$ such that $p(a)=\hat a$. 
Define $\rho= \hat \rho \circ p: G \to \S$.

Remark \ref{rem}(2) completes the proof. 
\end{proof}

Genus-1 knots are plentiful, the simplest examples being the twist knots (e.g. the knots $5_2, 6_1$) and doubled knots (obtained from a knot and any push-off by joining with a clasp). We extend the collection of nonfibered knots with uncountable representation shifts by considering also any knot $k$ with group $G$ that maps homomorphically onto the group $\bar G$ of a nonfibered genus-1 knot $\bar k$. Examples of such knots $k$ include satellite knots with genus-1 pattern knot \cite{whitten}.

\begin{cor} Let $k$ be a knot. Assume that the group of $k$ maps onto the group of a nonfibered knot $\bar k$ of genus 1.  Then $k$ is nonfibered and the conclusions of Proposition \ref{summary} hold. 

\end{cor}

\begin{proof} Assume that $h: G \to \bar G$ is an epimorphism, where $G$, $\bar G$ are the groups of  $k$, $\bar k$, respectively.   Let $K$, $\bar K$ denote the respective commutator subgroups, and $x$, $\bar x$ the meridianal generators of $k$, $\bar k$. 

Since $h(K)= \bar K$ and $\bar  K$ is not finitely generated, we see at once that $K$ is not finitely generated. Hence $k$ is nonfibered.

If $h(x) = \bar x$, then for any finite group $\S$, the representation shift $\Phi_\S(\bar\G)$ corresponding to $\bar k$ is a subshift of the representation shift $\Phi_\S(\G)$ corresponding to $k$; that is, ${\rm Hom}(\bar K, \S)$ is a subspace of ${\rm Hom}(K, \S)$ with the shift map $\sx$ restricting to $\sigma_{\bar x}$. The epimorphism  $h$ induces an embedding: $h^* \rho = \rho \circ h$. 
It follows that the topological entropy $h_\S(\G)$ is no less than $h_\S(\bar\G)$. By theorem \ref{main}, $h_\S(\bar \G)>0$ for some finite group $\S$. 
Hence for such a group, $h_\S(\G)$ is also positive. 

If $h(x) \ne \bar x$, then there exists $a \in K$ such that $h(a x) = \bar x^{\epsilon}$, where $\epsilon= \pm 1$. We may assume without loss of generality that $\epsilon=1$. 
In this case, we replace $x$ by $a x$. Of course the augmented group system $\G$ and associated representation shifts $\Phi_\S(\G)$ change. However, by a result of \cite{swIJM}, the topological entropy of the representation shift remains unchanged. As in the case in which $h(x)=\bar x$, there exists a finite group $\S$ such that $h_\S(\G)>0$. \end{proof}

 \bigskip

\noindent {\sl Address for both authors:} Department of Mathematics and  Statistics, ILB 325, University of South Alabama, Mobile AL  36688 USA \medskip

\noindent {\sl E-mail:} silver@jaguar1.usouthal.edu; swilliam@jaguar1.usouthal.edu

\end{document}